\documentclass[12pt]{amsart}

\usepackage{amscd}
\usepackage[all]{xy}
\usepackage{hyperref}
\usepackage{gensymb}

\newcommand{\PP}{\mathbb P}
\newcommand{\AAA}{\mathbb A}

\def\cB{{\mathcal B}}
\def\cC{{\mathcal C}}

\def\cO{{\mathcal O}}

\newtheorem{theorem}{Theorem}[section]

\theoremstyle{definition}
\newtheorem{definition}[theorem]{Definition}
\newtheorem{example}[theorem]{Example}

\theoremstyle{remark}

\numberwithin{equation}{section}

\begin{document}

\title[A Coarse Stratification of the Monster Tower]{A Coarse Stratification \\ of the Monster Tower}

\author[A.\ Castro]{Alex Castro}
\address{Lab49, 30 St. Mary Axe,
London EC3A 8EP,
UK}
\email{alex.decastro@lab49.com}

\author[S.\ J.\ Colley]{Susan Jane Colley}
\address{Department of Mathematics,
Oberlin College, Oberlin, Ohio 44074, USA}
\email{sjcolley@math.oberlin.edu}

\author[G.\ Kennedy]{Gary Kennedy}
\address{Ohio State University at Mansfield, 1760 University Drive,
Mansfield, Ohio 44906, USA}
\email{kennedy@math.ohio-state.edu}
\thanks{This work was partially supported by a grant from the Simons Foundation (\#318310 to Gary Kennedy).}

\author[C.\ Shanbrom]{Corey Shanbrom}
\address{California State University, Sacramento, 6000 J St., Sacramento, CA 95819, USA}
\email{corey.shanbrom@csus.edu}

\begin{abstract}
The monster tower is a tower of spaces over a specified base;
each space in the tower is a parameter space for curvilinear data up to a specified order.
We describe and analyze a natural stratification of these spaces.
\end{abstract}

\date{\today}
\maketitle

%%%%%%%%%%%%%%%%%%%%%%%%%%%%%%%%

\par
The monster tower, also known in the algebro-geometric literature as the Semple tower,
is a tower of smooth spaces (varieties) over a specified smooth base $M$.
Each space $M(k)$ in the tower, called a \emph{monster space},
is a parameter space for curvilinear data up to order $k$ on $M$. 
We will describe a coarse stratification of each monster space, 
with each stratum corresponding to a code word created 
out of a certain alphabet according to rules that we will specify.  
These strata parametrize curvilinear data ``of the same type.'' 
The monster space can be regarded as an especially nice compactification
of the parameter space for curvilinear data of nonsingular curves on $M$
(as explained in Section \ref{baby}), with the added points representing
the data of singular curves; in our stratification the nonsingular data points
will form a single open dense stratum.
\par
Versions of this coarse stratification have been observed by virtually everyone who has studied the monster construction. Here we develop the theory in full generality, beginning with a base space of arbitrary dimension and at all levels. A finer stratification would result from a thorough analysis of the orbits under the action of a suitable group acting on the base (whose action can be lifted to the tower) or, working locally, of the pseudogroup of local diffeomorphisms at a selected point. Results such as those in Section 5.7 of \cite{MZ1} show that one can expect there to be infinitely many strata, i.e., that there are moduli. We seem to be very far, however, from a full understanding of where and why moduli occur.
\par
In Section \ref{basics}
we recall the construction of the monster tower and selectively review prior literature.
A brief Section \ref{baby} explains curvilinear data and baby monsters.
Section \ref{codewords} introduces code words for labelling the strata,
which are explained in Section \ref{strata} via their closures,
 called \emph{intersection loci}.
Our main Theorem \ref{stratatheorem} gives an explicit description of these loci.
To prove the theorem we will use coordinates on charts, as
explained in Section \ref{charts}; after this the proof
in Section \ref{maintheoremproof} is nearly immediate.
Finally in Section \ref{counting} we count the strata.

%%%%%%%%%%%%
\section{The monster tower} \label{basics}

Suppose that $M$ is a smooth manifold, complex manifold, or nonsingular algebraic variety over
a field of characteristic zero; denote its dimension by $m$. The \emph{monster tower} is a sequence
\begin{equation} \label{towersequence}
\dots \to M(k) \to M(k-1) \to \dots \to M(2) \to M(1) \to M(0)=M
\end{equation}
in which each $M(k) \to M(k-1)$ is a fiber bundle with fiber
$\PP^{m-1}$.
To define it, we begin with a more general construction.
\par
The general construction begins 
with a pair $(X, \cB)$, where $X$ is again a smooth manifold, complex manifold, or variety,
and 
$\cB$ is a rank $b$ subbundle of its tangent bundle $TX$.  
Let $\widetilde{X}=\PP\cB$ and let $\pi : \widetilde{X} \to X$ be the projection. 
Then the tautological line bundle $\cO_{\cB}(-1)$ on $\widetilde X$ is 
a subbundle of $\pi^*TX$.  
Let
$$
d\pi : T{\widetilde{X}} \to \pi^*TX
$$ 
denote the derivative map of $\pi$.  
The pullback of $\cO_{\cB}(-1)$ by $d\pi$ is a subbundle of $T{\widetilde{X}}$, which we denote by $\widetilde \cB$.
In other words, a point $P$ of $\widetilde{X}$ represents a tangent direction at a point
of $X$, and $\widetilde \cB$ is the subbundle of $T{\widetilde{X}}$ whose
fiber at $P$ consists of vectors mapping (via the derivative of projection)
to vectors in that direction;
we call them \emph{focal vectors}.
Note that the relative tangent bundle  $T({\widetilde X / X})$ is a subbundle of
$\widetilde \cB$; its fiber consists of vectors mapping to zero,
sometimes called \emph{vertical vectors}.
By construction, $\widetilde \cB$ is a subbundle of $T{\widetilde{X}}$
and again its rank is $b$. Thus we can iterate this construction to obtain
a tower of fibrations.
\par 
We will eventually apply this construction in several situations.
To construct the monster tower, we simply apply it
to the pair $(M,TM)$ and then iterate. We denote the resulting spaces as in 
(\ref{towersequence});
$M(k)$ is called the \emph{monster space at level $k$}
or simply the \emph{$k$th monster}.
The bundle constructed at step $k$ of the construction is called the $k$th \emph{focal bundle}
and denoted $\Delta_k$; it is a subbundle of the tangent bundle $TM(k)$.
For $k \geq 2$,
the projectivization of the relative tangent bundle $\PP T(M(k)/ M(k-1))$
gives us a divisor on $M(k)$, which we call the 
\emph{divisor at infinity} and denote by $I_k$.
The pullback of $I_k$ to any higher level is again a divisor,
and for simplicity of notation we will again just denote it by $I_k$.
\par
The earliest instance of the construction seems to be Gherardelli's paper \cite{Gher}.
The tower was explained by Semple in \cite{Semple1, Semple2},
and in the algebro-geometric literature it bears his name.
Two of the present authors used it to study problems of enumerative geometry
in \cite{CK1, CK2, CK3}, and it was treated in greater generality in \cite{L-J}.
Demailly used it to study positivity questions for hyperbolic varieties in \cite{Dem}.
\par
Working independently, a group of differential geometers studied the same construction using different techniques, language, and motivation.  
In this strand of literature the general construction is called \emph{Cartan prolongation},
and the resulting tower is called the monster tower.
The tower with base $\mathbb R^2$ was discovered by Montgomery and Zhitomirskii \cite{MZ2} in their study of singular Goursat distributions.  They discovered connections with singular plane curves as well as the control-theoretic problem of a car pulling many trailers.  Their detailed study of the monster tower with base $\mathbb R^2$ appears as the monograph \cite{MZ1}.  Subsequent investigations of the monster tower with base $\mathbb R^2, \mathbb R^3$, and $\mathbb C^3$ appear in \cite{castro2, castro3, CHS, CS}.  These efforts were generally aimed at understanding the action of the diffeomorphism group on the tower: constructing invariants and counting and classifying orbits.
\par
The equivalence of the two towers was first noticed by Castro in 2010. The present contribution represents the first collaboration between the two groups, and an effort to improve and standardize language and notation, most notably in the coding system presented in Section \ref{codewords}.

\section{Curvilinear data and baby monsters} \label{baby}
Here we informally recall the concept of curvilinear data; for further information see the
works cited at the end of the previous section.
\par
Again let $M$ be a smooth manifold, complex manifold, or nonsingular algebraic variety
of dimension $m$. 
Suppose we have two smooth curves $C_1$ and $C_2$ passing through
a point, and that we have a system of local coordinates 
$x_1, \dots, x_m$
based there.
For each curve, assume that the restriction of the differential 
$dx_1$ does not vanish at the point.
We say that the curves \emph{have the same curvilinear data}
up to order $k$
at the point if
the values of all derivatives
$d^{j}x_i/dx_1{}^{j}$
agree up to order $k$.
These are $(m-1)k$ conditions, and
one can check that they are independent of local coordinates.
\par
Thus a 
\emph{nonsingular curvilinear datum} is a point
in a manifold or smooth variety of dimension $m+(m-1)k$,
and there is a tower 
of such manifolds over the base
$M$, with fiber a projective space $\PP^{m-1}$
at the first level, and then affine space fibers $\AAA^{m-1}$
thereafter.
At the first level the manifold is exactly $M(1)$,
the projectivization of the tangent bundle $TM$.
To each point of the curve $C$ in $M$ we can associate
the point of $M(1)$ recording its tangent direction,
and in this way we obtain a copy of $C$ inside $M(1)$,
called the \emph{lift} of $C$ and denoted by $C(1)$; the process of passing from
$C$ to $C(1)$ is again called \emph{prolongation}.
\par
We can lift again to obtain a curve $C(2)$ inside the projectivization
of the tangent bundle $TM(2)$. But observe
that a tangent vector at a point of $C(1)$
is a focal vector; thus in fact $C(2)$ lies entirely inside the smaller space
$M(2)$. Proceeding similarly, we obtain a copy $C(k)$
of $C$ inside the monster space $M(k)$,
called its \emph{$k$th lift} or \emph{$k$th prolongation}.
At every stage these focal vectors are non-vertical vectors;
thus $C(k)$ avoids the divisor at infinity $I(k)$, as well
as all pullbacks of prior divisors at infinity $C(2)$, \dots, $C(k-1)$.
Conversely, given any point of $M(k)$ away from all divisors at infinity,
there is a recipe for finding a curve $\widetilde{C}$ passing through this point and
then, by a process of integration, producing a curve $C$ on $M$
for which $C(k)=\widetilde{C}$. Thus the monster spaces $M(k)$
are spaces which naturally compactify the spaces of nonsingular curvilinear
data. As the reader undoubtedly suspects, a point on a divisor
at infinity represents the data of some sort of singular curve;
the process of lifting is essentially repeated Nash blowup, performed
simultaneously at every point of the curve.
The added points are also called \emph{curvilinear data},
dropping the modifier ``nonsingular.''
This is explained more carefully in the cited literature. 
\par
Inside the monster space $M(k)$, 
suppose that there is a submanifold $X$ for which the intersection of its tangent bundle $TX$
with the focal bundle $\Delta_k$ (inside the tangent bundle of $M(k)$ restricted to $X$)
has constant rank. Then we can apply the monster prolongation construction to the 
pair $(X,TX\cap \Delta_k)$; we call the resulting tower the \emph{baby monster tower}
associated to $X$.
\par
Here are three basic examples:
\begin{enumerate}
\item
If $C$ is a smooth curve in $M$ and $TC$ is its tangent bundle,
then applying the baby monster construction to the pair $(C,TC)$
produces a tower of copies of $C$, namely the lifts $C(k)$ in $M(k)$,
as just explained.
\item
If $X$ is the fiber over a point of $M(k-1)$, then
$TX$ is already a subbundle of $\Delta_k$.
The resulting tower is what Castro and Montgomery (\cite{castro2})
have also called a baby monster.
(Thus our terminology generalizes theirs.)
\item
If $X$ is the divisor at infinity $I_k$
on $M(k)$,
then the intersection of $TX$ and $\Delta_k$
is transverse; thus $TX\cap \Delta_k$ has rank $m-1$. This tower is
described by Lejeune-Jalabert \cite{L-J} on page 1287.
Denote the spaces in this tower by $I_k[n]$,
beginning with $I_k[0]=I_k$;
thus $I_k[n]$ is a subspace of $M(k+n)$
of codimension $n+1$.
We use the same notation to denote the full inverse image
of $I_k[n]$ in any higher monster space in the tower.
\end{enumerate}
\par

%%%%%%%%%%%%
\section{Code words} \label{codewords}
\par
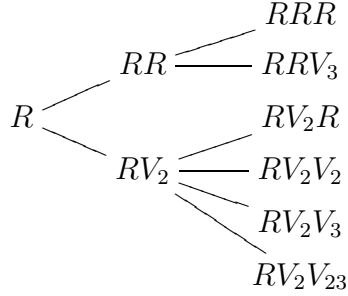
\begin{figure}
\[
\xymatrix@R=0.25em{ 
              &       & RRR \ar@{-}[dl]\\
              &   RR \ar@{-}[dl]                   & RRV_3 \ar@{-}[l] \\
R            &                      & RV_2R \ar@{-}[dl] \\ 
              &           RV_2 \ar@{-}[ul]            &  RV_2V_2 \ar@{-}[l]\\
              &                      & RV_2V_3\ar@{-}[ul] \\ 
              &                      & RV_2V_{23} \ar@{-}[uul] \\ 
               }
\]
\caption{Code words up to length 3}
 \label{wordstolength3}
\end{figure}

We now introduce the code words to be used for labeling strata;
the strata themselves will be described later.
The alphabet for our code consists of all symbols $V_A$,
where $A$ is a finite subset of the integers strictly greater than 1.
To avoid a cumbersome notation, we write, e.g., $abc$ or  $a,b,c$ rather than $\{a,b,c\}$,
and always arrange the elements of $A$ in increasing order.
Although this is an infinite alphabet, the rules for creating a valid word
will imply that there are only finitely many words of each specified length.
To be consistent with prior usage, we will use $R$ in place of $V_{\emptyset}$.
(The symbols $R$ and $V$ have been chosen to suggest ``regular'' and ``vertical.'')
The rules for creating a code word are as follows:
\begin{enumerate}
\item
The first symbol must be $R$.
\item
Immediately following the symbol $V_A$, one may put any symbol $V_B$, where
either $B$ is a subset of $A$, or $B$ is a subset of $A\cup\{j\}$, with $j$
being the position of the symbol.
\item \label{vlength}
The cardinality of  $A$ is less than $m$.
\end{enumerate}
Note that $j$ cannot appear in a subscript prior to position $j$.
Also note that rule (\ref{vlength}) is the only rule to use the specified value for $m$.
The diagram in  
Figure~\ref{wordstolength3}
shows the code words of lengths 1, 2, and 3,
assuming that $m \geq 3$.
If $m>3$, then there are twenty-four valid code words of length four.
If $m=3$, however, then
$RV_{2}V_{23}V_{234}$ is not a valid code word,
since it violates rule (\ref{vlength}).
\par

Given a code word $V_{A_1}V_{A_2}V_{A_3}\cdots V_{A_k}$, for each $j=2,3,\dots,k$,  let $n_j$ denote the number of times that $j$ appears as a subscript.
(If it appears, then its first appearance is in position $j$ and its last appearance is in position $j+n_j-1$.)
For example, for the code word $RV_2V_{23}V_{23}V_{25}V_5V_5V_5$ we have 
$n_2=4$, $n_3=2$, $n_4=0$, and $n_5=4$. 
We note that $n_j \leq k+1-j$; the base dimension $m$ adds further restrictions.
From a valid set of values for $n_2$ through $n_k$
we can recover the code word. 
\par
In \cite{castro2} and \cite{CHS}, an alternative coding system was developed for the case $m=3$, 
called \emph{RVT coding}. In this system one works with the 7-symbol alphabet
$\{R,V,T_1,T_2,L_1,L_2,L_3\}$, which corresponds to the seven cases considered by Semple on
page 151 of \cite{Semple1}. The rules for
creating a code word are as follows:
\begin{itemize}
\item
The first symbol must be $R$.
\item
The next symbol after $R$ can be either $R$ or $V$.
\item
The next symbol after $V$ or $T_1$ can be $R$, $V$, $T_1$, or $L_1$.
\item
The next symbol after $T_2$ can be $R$, $V$, $T_2$, or $L_3$.
\item
The next symbol after $L_1$, $L_2$, or $L_3$ can be any symbol.
\end{itemize}
\par
The symbol $R$ plays the same role in both the
$RVT$ code and our new code.
In the correspondence between the two codes, $V$, $T_1$, and $T_2$
always correspond to singly-subscripted symbols, but
the precise correspondence is conditioned by prior symbols of the code word;
similarly $L_1$, $L_2$, and $L_3$ correspond to doubly-subscripted symbols.
The precise correspondence is explained by these tables,
in which the first row refers to the $RVT$ code and the second row to our new code.

\bigskip
%%Table1
\begin{center}
\begin{tabular}{ r |c |c |c |c |c |c |c |}
\cline{2-3}
$k$th symbol, immediately after $R$ 
	& $R$ & $V$   \\ \cline{2-3} 
$k$th symbol, immediately after $R$ 
	& $R$ & $V_k$ \\ \cline{2-3}
\end{tabular} 
\bigskip

%%Table2
\begin{tabular}{r |c |c |c |c |c |c |c |}
\cline{2-5}
$k$th symbol, after $V$ or $T_1$ 
	& $R$ & $V$ & $T_1$ & $L_1$    \\ \cline{2-5}
$k$th symbol, after $V_{j}$ 
	& $R$ & $V_k$ & $V_j$ & $V_{jk}$ \\ \cline{2-5} 
\end{tabular} 
\bigskip

%%Table3
\begin{tabular}{r |c |c |c |c |c |c |c |}
\cline{2-5} 
$k$th symbol, after $T_2$ 
	& $R$ & $V$ &  $T_2$ & $L_3$   \\ \cline{2-5} 
$k$th symbol, after $V_{j}$ 
	& $R$ & $V_k$ & $V_j$ & $V_{jk}$ \\ \cline{2-5}
\end{tabular} 
\bigskip

%%Table4
\begin{tabular}{r |c |c |c |c |c |c |c |}
\cline{2-8}
$k$th symbol, after $L_1$, $L_2$, or $L_3$ 
	& $R$ & $V$ & $T_1$ & $T_2$ & $L_1$ & $L_2$ & $L_3$   \\ \cline{2-8} 
$k$th symbol, after $V_{ij}$ 
	& $R$ & $V_k$ & $V_j$ & $V_i$ & $V_{jk}$ & $V_{ij}$ & $V_{ik}$ \\ \cline{2-8} 
\end{tabular} 
\end{center}
\medskip
\noindent
It is straightforward to use these tables to translate an $RVT$ code word to a new
code word. Here is one example:
$$
R V T_1 V V T_1 T_1 T_1 L_1 \longleftrightarrow R V_2 V_2 V_4 V_ 5 V_5 V_5 V_5 V_{59}.
$$
In the other direction, when translating the symbol of a new code word
immediately following some $V_j$, one one must be careful to choose the correct table.
For example, although $RV_2V_{23}V_2V_2$ and $RV_2V_2V_2V_2$ end
with the same pair of symbols, the final symbols of their corresponding $RVT$
code words differ:
$$
\begin{aligned}
RV_2V_{23}V_2V_2 &\longleftrightarrow RVL_1T_2T_2, \\
RV_2V_2V_2V_2 &\longleftrightarrow RVT_1T_1T_1.
\end{aligned}
$$
 
%%%%%%%%%%%%
\section{Intersection loci and strata} \label{strata}
We now give a recipe for converting a code word 
$V_{A_1}V_{A_2}V_{A_3}\cdots V_{A_k}$ to a description
of an intersection locus. Recall that for each number $j=2,3,\dots,k$,
we let $n_j$ denote the number of times that $j$ appears as a subscript.
If $j$ never appears, then set $n_j=-1$.
In order to state Definition~\ref{intlocdef} cleanly, we adopt the following convention:
$I_j[-1]$ is interpreted as the entire monster space $M(k)$.

\begin{definition}
\label{intlocdef}
On $M(k)$, the intersection locus corresponding to the
code word $W=V_{A_1}V_{A_2}V_{A_3}\cdots V_{A_k}$ is
\begin{equation} \label{intloceq}
I_W:=\bigcap_{j=2}^k I_j[n_j-1].
\end{equation}
\end{definition}
\begin{theorem} \label{stratatheorem}
For each code word,
the intersection (\ref{intloceq}) is transverse and nonempty.
The codimension of the intersection locus is 
the sum of cardinalities
$|A_1|+|A_2|+\cdots+|A_k|$, equivalently $n_1+n_2+\cdots+n_k$.
It contains the intersection locus 
corresponding to the code word $W'$ if and only if each $n'_j \geq n_j$.
\end{theorem}

In Section \ref{maintheoremproof}
we will prove Theorem~\ref{stratatheorem} by an explicit calculation, using 
coordinates to be introduced in Section~\ref{charts}.
\par
\begin{example}
For the code word consisting entirely of $R$'s,
the intersection locus is $M(k)$ itself.
\end{example}
\begin{example}
Assume $m \geq 3$. 
On the monster $M(8)$,
the codimension-7 intersection locus corresponding to the code word 
$RV_2V_3V_{34}V_{35}V_3RR$
is the transverse intersection
$$
I_2[0] \cap I_3[3] \cap I_4[0] \cap I_5[0],
$$
where $I_j[0]$ means the full inverse image of $I_j$ and $I_3[3]$ denotes the inverse image of the third space in the baby monster tower over $I_3$.  
\end{example}
\par
Using Theorem~\ref{stratatheorem}, we obtain a natural stratification of the monster space $M(k)$
as follows: from each intersection locus $I_W$ excise
all those intersection loci $I_{W'}$ that it contains.
In fact, it suffices just to excise each locus whose code word $W'$
is obtained from $W$ by increasing a single $n_j$ by one
(ignoring those $W'$ that are not valid code words).

%%%%%%%%%%%%
\section{Coordinate charts} \label{charts}
Here we describe a natural system of coordinate 
charts on the monster spaces.
The notation is basically that of Lejeune-Jalabert \cite{L-J},
except for the conveniently redundant names.  Note that a similar system was developed from the differential geometry side in \cite{GKR} (see \cite{MZ1} for more details).
\par
We begin with an open set $U$ on $M$
with coordinates $x_1, \dots, x_m$ for which, at each point, the differentials
$dx_1, \dots, dx_m$ form a basis of the cotangent space.
Then over $U$ the monster space $M(1)$ is 
isomorphic to $U \times \PP^{m-1}$, and
it is covered by $m$ charts $\cC(1), \dots, \cC(m)$, 
each isomorphic to $U \times \AAA^{m-1}$. 
The chart $\cC(p)$ represents one-dimensional quotients of the cotangent bundle
on which $dx_p$ does not vanish, and the coordinates for the second factor
are defined by 
$$
x_q(p):=\frac{dx_q}{dx_p}
$$
for $q \neq p$.
For convenience we also
define
$$
x_p(p):=x_p
$$
(i.e., we give a new name to the pullback of this coordinate function).
Note that at each point in this chart the differentials
$dx_1(p), \dots, dx_m(p)$ form a basis for the 
dual focal bundle $\Delta^{\vee}_1$.
For this reason the coordinates $x_q(p)$ are called 
\emph{active coordinates};
the active coordinate $x_p(p)$ is also called the \emph{retained coordinate}.
\par
This is the beginning of a recursive construction of
a system of
$m^k$ charts $\cC(p_1p_2\dots p_k)$ on $M(k)$.
Each chart is isomorphic to
$$
\cC(p_1p_2\dots p_{k-1}) \times  \AAA^{m-1},
$$
and there are $m$ \emph{active coordinates}
$x_q(p_1p_2\dots p_{k})$, where $1 \leq q \leq m$.
The \emph{retained (active) coordinate} is 
$$
x_{p_k}(p_1p_2\dots p_{k}):=x_{p_k}(p_1p_2\dots p_{k-1})
$$
(the pullback of a coordinate from below, given a new name);
the others are defined by
$$
x_{q}(p_1p_2\dots p_{k}):=\frac{dx_{q}(p_1p_2\dots p_{k-1})}{dx_{p_k}(p_1p_2\dots p_{k-1})}
$$
and serve as affine coordinates for the second
factor $\AAA^{m-1}$.
The chart represents 
one-dimensional quotients of $\Delta^{\vee}_{k-1}$
on which the differential of the retained coordinate does not vanish.
One verifies that at each point the differentials of the active coordinates
form a basis for $\Delta^{\vee}_{k}$.
\par
Thus for each chart we have shown how to systematically construct $(k+1)m$ coordinate names
$$
\begin{aligned}
x_1&, \dots, x_m, \\
x_1&(p_1), \dots, x_m(p_1), \\
x_1&(p_1p_2), \dots, x_m(p_1p_2), \\
&\vdots \\
x_1&(p_1p_2\dots p_k), \dots, x_m(p_1p_2 \dots p_k).
\end{aligned}
$$
Because of redundant names, there are in fact $m+k(m-1)$
distinct coordinates.
For a coordinate with redundant names,
there is a \emph{shortest name},
characterized by the fact that its subscript
differs from the final symbol appearing in parentheses. 
(The names $x_1, \dots, x_m$ are always shortest names.)

%%%%%%%%%%

\begin{example} \label{chartexample}
Assuming that $m=3$, in chart $\cC(32123)$
on $M(5)$ we have the following coordinate names:
\begin{align*}
&x_1 && x_2 && x_3 \\
&x_1(3) && x_2(3) && \boxed{x_3(3)} \\
&x_1(32) && \boxed{x_2(32)} && x_3(32) \\
&\boxed{x_1(321)} && x_2(321) && x_3(321) \\
&x_1(3212) && \boxed{x_2(3212)} && x_3(3212) \\
&x_1(32123) && x_2(32123) &&  \boxed{x_3(32123)}
\end{align*}
All names are shortest names except those in the boxes.
Here $x_3(3212)$ is the retained active coordinate.

\end{example}

%%%%%%%%%%%

\par
We now identify, in each chart $\cC(p_1p_2\dots p_k)$, the equations
for the loci appearing 
in Theorem~\ref{stratatheorem}.
The divisor at infinity $I_j$ first 
appears on $M(j)$.
It represents one-dimensional quotients of $\Delta^{\vee}_{j-1}$
in which the differential of every coordinate pulled back
from $M(j-2)$ vanishes. 
There are two possibilities. If $p_{j-1} = p_j$, 
then the retained coordinate is pulled back from $M(j-2)$; 
thus $I_j$ does not meet the chart $\cC(p_1p_2\dots p_j)$.
If $p_{j-1} \ne p_j$, then we claim that the differential of every coordinate 
pulled back from $M(j-2)$ vanishes if and only if 
\begin{equation}\label{divinfequation}
x_{p_{j-1}}(p_1p_2\dots p_{j})=0.
\end{equation}
Indeed, this equation is satisfied if and only if the differential
of the previous retained coordinate vanishes:
$$
dx_{p_{j-1}}(p_1p_2\dots p_{j-1})=0,
$$
and in $\Delta^{\vee}_{j-1}$ the differentials of all inactive
coordinates are multiples of this differential.
\par
Thus on the divisor at infinity $I_j$ the differential 
$dx_{p_{j-1}}(p_1p_2\dots p_{j})$ vanishes. 
When we prolong, there are again two cases.
If $p_{j-1} = p_{j+1}$, then again $I_j[1]$ does not meet the chart 
$\cC(p_1p_2\dots p_{j+1})$.
If $p_{j-1} \ne p_{j+1}$, then when we prolong we have
\begin{equation}\label{prolongdivinfequation}
x_{p_{j-1}}(p_1p_2\dots p_{j}p_{j+1})=0,
\end{equation}
so that $I_j[1]$ is defined by both (\ref{divinfequation})
and (\ref{prolongdivinfequation}).
Continuing in this manner, we see that $I_j[n_j-1]$
is defined by the vanishing of $n_j$ coordinates
$$
x_{p_{j-1}}(p_1p_2\dots p_{j}), \quad
x_{p_{j-1}}(p_1p_2\dots p_{j+1}), \quad
\dots, \quad
x_{p_{j-1}}(p_1p_2\dots p_{j+n_j-1}),
$$
always assuming that the symbol $p_{j-1}$ is not repeated
later in the word; 
if it is repeated, then $I_j[n_j-1]$ does not meet the chart. 
\par
\begin{example}
Suppose $m \ge 3$ and consider the code word $RV_2V_2V_{24}R$, for which $n_2 = 3$, $n_3 = 0$, $n_4 = 1$, and $n_5 = 0$.  It represents the stratum on $M(5)$ whose closure is the intersection locus $I_2[2] \cap I_4$; to obtain the stratum one needs to excise these smaller intersection loci:
\[   I_2[3] \cap I_4, \quad I_2[2] \cap I_3 \cap I_4, \quad I_2[2] \cap I_4[1],  \quad I_2[2] \cap I_4 \cap I_5.    \]
This stratum has codimension 4.

\par
Now assume $m=3$.  Here are the equations of $I_2[2] \cap I_4$ in the chart $\cC(32123)$, using the coordinate names shown in Example~\ref{chartexample}:
\begin{align*}
x_3(32)&=0\\
x_3(321)&=0\\
x_3(3212)&=0\\
x_1(3212)&=0.
\end{align*}
The first three equations define $I_2[2]$, while the fourth defines $I_4$.  The following table displays the additional equation needed to define each of the excised intersection loci:
\medskip
\begin{center}
\begin{tabular}{lll}
$I_2[3] \cap I_4$ & & does not meet the chart \\[4pt]
$I_2[2] \cap I_3 \cap I_4$ & & $x_2(321) = 0$ \\[4pt]
$I_2[2] \cap I_4[1]$ & & $x_1(32123) = 0$ \\[4pt]
$I_2[2] \cap I_4 \cap I_5$ & &  $x_2(32123) = 0$.
\end{tabular}
\end{center}

\end{example}

%%%%%%%%%%%%
\section{Proof of the main theorem} \label{maintheoremproof}
We now prove Theorem~\ref{stratatheorem}. 
If the intersection locus 
\begin{equation} \label{intersections}
\bigcap_j I_j[n_j-1]
\end{equation}
meets a chart $\cC(p_1p_2\dots p_k)$,
then (as we have explained in Section \ref{charts})
it is defined by the vanishing of $\sum_{j=1}^k n_j$ coordinates, noting that the names we
have given these coordinates are all shortest names and so there are no repetitions.
Thus the intersection is transverse and the codimension of the intersection locus is as claimed. 
These explicit equations also make clear the claim about when one intersection locus contains another.
The stratum for the code word $V_{A_1}V_{A_2}V_{A_3}\cdots V_{A_k}$ is obtained from 
(\ref{intersections}) by removing all the smaller intersection loci.
\par
Finally we argue that each intersection locus, and hence each stratum, is nonempty.
Indeed, suppose that $V_{A_1}V_{A_2}V_{A_3}\cdots V_{A_k}$ is a valid code word.
Then the chart $\cC(p_1p_2\dots p_k)$
will meet the intersection locus (\ref{intersections})
if and only if each subscript $p_j$ avoids a certain subset of the 
previous subscripts, and this subset has the cardinality of $A_j$.
Since this cardinality is less than $m$, one can always choose such a $p_j$.

%%%%%%%%%%%%

\section{Counting code words} \label{counting}
In this section we explain how to count the number of code words.
As before, fix the dimension of the base space to be $m$.
Let $N(k,r)$ be the number of code words of length $k$ 
in which the last subscript has length $r$.
Assume that $k \geq 2$ and that $1 \leq r \leq k-1$.
Then
$$
N(k,r) = \sum_{i=r-1}^{m-1} \binom{i+1}{r} N(k-1,i).
$$
Indeed, to obtain such a code word, take any code word of length $k-1$
in which the last subscript $A$ has length $i \geq r-1$;
then create the new subscript $B$ 
by choosing any $r$ symbols from the set $A\cup\{k\}$.
Clearly
$$
N(k,0) = \sum_{i=0}^{m-1} N(k-1,i) = \text{number of code words of length $k-1$}.
$$
To use these recursive formulas, begin with $N(1,0)=1$ and $N(1,r)=0$ for $r>0$.
\par
We claim that if $k \leq m$, then 
$N(k,r)$ equals the unsigned Stirling number of the first kind
$c(k,r+1)$.
One way to establish this is to remark that these Stirling numbers satisfy the same recurrence
just established for the number of code words.
\par
Alternatively, one can establish a 
bijection between the set of code words of length $k$ 
and the set of all rooted trees on the vertex set $\{0, 1, \dots, k\}$
in which $0$ is the root and the labels increase as one moves
away from the root (i.e., \emph{increasing trees}).
Given a code word, recall that $n_j$ denotes the number
of times that $j$ appears as a subscript.
If the code word has length $k$, 
create a tree with vertices labeled by the integers $1, 2, \dots, k+1$,
with root at $k+1$. 
For each integer $j=2, \dots, k$, draw
an edge connecting $j-1$ and $j+n_j$;
also draw an edge connecting $k$ and $k+1$.
Note that the degree of the root
is $r+1$,
and that 
labels increase as one moves toward the root.
Replacing each label $j$ by $k+1-j$,
one obtains an increasing tree.
The process clearly can be reversed;
thus we have the desired bijection.
By Proposition 1.5.5(b) of \cite{Stanley},
$N(k,r)=c(k,r+1)$.

\bibliographystyle{amsplain}

 \end{document}